\documentstyle[11pt]{article}
\font\cyr=wncyr10
\def\Sh{\hbox{\cyr Sh}}
\begin{document}
\title{ {Mordell-Weil Groups and Selmer Groups of
 \\ Two Types of Elliptic Curves }
     \thanks{MR(1991) Subject Classification:
      11R58; 11R27; 14H05; 11A55}
         \thanks{Project Supported by the NNSFC (No. 19771052)}
         }
\author{\mbox{}\vspace{0.3cm}
   {{\Large Q}IU {\Large D}erong  $\;\;
               and\;\; $ {\Large Z}HANG {\Large X}ianke }\\
 \small{Tsinghua University,
     Department of Mathematical Sciences,
            Beijing 100084, P. R. China  } }
\date{}
\maketitle
\parindent 24pt
\baselineskip 18pt
\parskip 2pt

 \begin{center}
 \begin{minipage}{14cm}
\baselineskip 18pt
\quad  Consider elliptic curves
$\ E=E_\sigma : \; y^2 = x (x+\sigma p)
     (x+\sigma q),\; $  where$\ \sigma =\pm 1,\; $
    $p$ and $ q$ are prime numbers with $p+2=q$.
 (1) The Selmer groups $\ S^{(2)}(E/{{\bf Q}}),\;
 S^{(\varphi)}(E/{{\bf Q})}$,  and
$\ S^{(\widehat{\varphi})}(E/{{\bf Q})}\ $
  are explicitly determined, e.g.,
    $\ S^{(2)}(E_{+1}/{{\bf Q}})= $
 $({{\bf Z}}/2{{\bf Z}})^2;\ $
 $ ({{\bf Z}}/2{{\bf Z}})^3;\ $ or
 $  ({{\bf Z}}/2{{\bf Z}})^4\ $
  when $p\equiv 5; \ 1 $ or $3;\ $ or
  $ 7 \ ({{\rm mod}} \ 8)$  respectively.
  (2) When $p\equiv 5  (\ 3 ,\ 5$
   for $\sigma =-1) \ ({{\rm mod}} \ 8),\ $
  it is proved that the Mordell- Weil group
  $\ E({{\bf Q})} \cong $
   $ {{\bf Z}}/2{{\bf Z}}
  \oplus{{\bf Z}}/2{{\bf Z}}\ $  having rank $0,\ $ and
  Shafarevich-Tate group  $\Sh(E/{{\bf Q}} )[2]=0.\ $
   (3) In any case,  the sum of
  rank$E({{\bf Q})}$ and dimension of
   $\Sh(E/{{\bf Q}} )[2] $ is given, e.g., $0;\ 1;\ 2 $
 when $p\equiv 5; \ 1 $ or $3;\ 7 \ ({{\rm mod}} \ 8)$
 for $\sigma =1$.
 (4) The Kodaira symbol, the torsion subgroup $E(K)_{tors}$
for any number field $K$ , etc. are also obtained.

{\bf{key words}}:
 elliptic curve, Selmer group, Mordell-Weil group,
    Shafarevich group
 \end{minipage}
\end{center}

 \par \vskip 0.66cm
\begin{center}
{\Large I.\ Introduction and Main Results}
\end{center}

 Let $p, q$  be two (twin) prime numbers and $q-p=2$ .
  Here we consider the elliptic curves
   $$\hskip 2cm E=E_\sigma :\; \; y^2 = x (x+\sigma p)
(x+\sigma q)\qquad (\sigma =\pm 1) \hskip 3cm (1.1)  $$
 We also denote $\ E=E_+$  or $E_-\ $   when $\sigma = +1$
 or $-1$.
 One of the interests to consider elliptic curves (1.1) is
 on twin primes.  Whether there are
infinitely many twin primes $p, q$ now is equivalent to
 whether there are infinitely many
isomorphic classes of such elliptic curves $E$ ,
since these elliptic curves   are not
isomorphic to each other for different $(p, q )$
as we will see later.   \par
Similarly to (1.1), some other special types of elliptic
curves  were studied by A. Bremner, J. Cassels, R. Strocker,
J. Top,  B. Buhler, B.  Gross and D.  Zagier
 (see [1-6]), $\; e.g.,\;  \; y^2=x(x^2+p), \; $
   $\; y^2=(x+p)(x^2+p^2) , \; $  and
   $\; y^2 = 4x^3 - 28x + 25 .\ $
   The first two elliptic curves have ranks
    $0,\ 1, \ $ or $\ 2. \ $ And the last elliptic curve has
   rank 3 and is famous in solving the Gauss conjecture.
\par
  For elliptic curves  $E$ in (1.1), it is easy to see
      the discriminate of it is
 $\ \Delta = \Delta(E) = 64p^2q^2.\  $
 the $j-$invariant is
 $\ j = j(E) = 64 (p^2+2q)^3/p^2q^2\ $
(Thus elliptic curves $E$  are not isomorphic to
 each other for different $(p, q )$). We
will  also show equation (1.1) is  global minimal.

 For any number field  $K$ and any elliptic curve $E$
  over $K$ ,  The Mordell-Weil theorem said
that the set $E(K)$ of  $K-$rational points of $E$
is a finitely generated abelian group (the Mordell-Weil group),
 so $\; \; E(K) \cong E(K)_{tors} \oplus {{\bf Z}}^r,\; $
 where $\ E(K)_{tors}\ $  is the torsion subgroup of $E(K)$ ,
 $\; r={{\rm rank}}E(K)\ $ is the rank of $E(K)$.
 We first consider $\ E(K)_{tors}\ $ for any number field $K$.

\par \vskip 0.36cm

{\bf Theorem 1}. Let $K$  be any number field. We have
the following results  on the torsion
subgroup $\ E(K)_{{tors}}\ $  for elliptic curve $E$
in (1.1) with $\ (p, q )\neq (3, 5),\
$  here $\ \wp\  $ is a prime ideal of $K$ over 3,
$\; e=e(\wp|3)\ $ and $\; f=f(\wp|3)\
$ are the ramification index and residue degree of
$\wp $ respectively. \par
 (1) If $\ e(\wp|3) = f(\wp|3)=1,\  $  then
  $$E(K)_{tors} \cong {{\bf Z}}/2 {{\bf Z}}
   \oplus {{\bf Z}}/2 {{\bf Z}}.$$
   \qquad (2) If  $\  f(\wp|3)=1 \ $  and $E$
   has an additive reduction
   (at any finite valuation of $K$ ), then
   $$E(K)_{tors} \cong {{\bf Z}}/2 {{\bf Z}}
  \oplus {{\bf Z}}/2 {{\bf Z}}
        \quad {\rm or} \quad {{\bf Z}}/2
 {{\bf Z}} \oplus {{\bf Z}}/6 {{\bf Z}}.$$
 \qquad (3) If   $\ f(\wp|3)=1\ $, then
  $$E(K)_{tors}/E(K)_3 \cong {{\bf Z}}/2 {{\bf Z}}
   \oplus {{\bf Z}}/2 {{\bf Z}},$$
where $E(K)_3$  denotes the 3-Sylow subgroup of $E(K)$ ,
 i.e. points of order a power of 3.
 \par     (4) If $E$  has an additive reduction
 (at any valuation),  then
 $\ \#E(K)_{tors} = 2^m \ $ or $\ 2^m\cdot 3 $ .
        \par \vskip 0.2cm

 We now introduce the Kodaria-Neron classification of the
special fibers ${{\cal C}}_\ell$  on the Neron
models of elliptic curves  $E/{{\bf Q}}$ .
 Let $\ell $  be a prime number,
 ${{\bf Q}}_\ell$  the $\ell-$adic
rationals (the completion of ${{\bf Q}}$ at the
 $\ell-$adic valuation $v_\ell$ ),
 $\; {{\bf F}_\ell = {{\bf Z}}/\ell {{\bf Z}}}$
 the finite field with $\ell $ elements. Let
   $\rho \ : \quad E({{\bf Q}}_\ell)
    \rightarrow \tilde{E}({{\bf F}}_\ell), \; \; $
   $P\mapsto \tilde{P}$
     be the reduction of $E$  modulo $\ell$ ,
      $\; \tilde{E}_{ns}({{\bf F}}_\ell)$
       the non-singular
${{\bf F}}_\ell-$points of
$\ \tilde{E}({{\bf F}}_\ell)$,
  $E_0({{\bf Q}}_\ell) =
  \rho^{-1}(\tilde{E}_{ns}({{\bf F}}_\ell)).\;$
 The Kodaira symbols
  ${I_0 ,\ I_1,\   I_n,\   II,\  III,\   IV,\ \cdots  }$
 are used to describe the type of the special
fiber ${{\cal C}}_\ell$  of the minimal
 Neron model of $E$  at $\ell$ ; $m_\ell$
denotes the number of irreducible components
 (ignoring multiplicities) on the special
fiber ${{\cal C}}_\ell$ .
 The conductor of $E/{{\bf Q}}$
 is defined by $N_E = \prod \ell ^{f_\ell}$ ,
where the product takes for $\ell $
 running over prime numbers,
 the local exponent $f_\ell $ is
defined by  $f_\ell  =0,\ 1,\  $ or
$\ 2+\delta_\ell$ according to $E$ having good,
multiplicative, or additive reduction at $\ell$,
 $\delta_\ell$ is concerned the
 action
of inertia group ( $\delta_\ell =0$ if
 $\ \ell \neq 2, 3 $  ). The index
 $\; c_\ell =| E({{\bf Q}}_\ell) / E_0({{\bf Q}}_\ell)| $
  is said to be  Tamagawa constant at $\ell$ .
 Using  results of Tate, Ogg, Kodaria-Neron, and Diamond
 (see [5-8]), we could obtain the following Theorem 2 and
Corollary 1(See [5-7] for notations).    \par \vskip 0.36cm

      {\bf Theorem 2}.  Let $E$  be an elliptic curve as
 in (1.1). Then we have the following Table.\par \vskip 0.18cm

\begin{tabular}{|c|c|c|c|c|}\hline
     Prime number $ \ell $
      & $ \;\ell=2 \; $ & $\; \ell=p\;  $ & $\; \ell=q\;  $
   & $\ell \neq2, p, q $\\
 \hline Kodaira Symbol
   & $\; III\; $ & $\; I_2\; $ & $\; I_2\; $ & $I_0 $\\
 \hline Number of Irreducible Component $ m_\ell\; $
   & $\; 2\; $ & $\; 2\; $ & $\; 2\; $ & $ 1  $ \\
 \hline Tamagawa Constant $ c_\ell $
   & $\; 2\; $ & $\; 2\; $ & $\; 2\; $ & $\; 1\; $ \\
 \hline $ E({{\bf Q}}_\ell) / E_0({{\bf Q}}_\ell)$
  & $\; \; {{\bf Z}}/2{{\bf Z}}\; \; $
   & $\; \; {{\bf Z}}/2{{\bf Z}}\; \; $
    & $\; \; {{\bf Z}}/2{{\bf Z}}\; \; $ & $ 0 $ \\
 \hline Local exponent $ f_\ell $  of conductor $ N_E $
   & $\; 5\; $ & $\; 1\; $ & $\; 1\; $ & $ 0 $ \\
  \hline
\end{tabular}

  \par \vskip 0cm

     {\bf Corollary 1}. The elliptic curves in (1.1) are
     {\bf modular} with conductor $N_E = 2^5pq,\;$  and their
 $L-$function $L(E/{{\bf Q},\ s}) $
  could be continually extended to the whole complex plane
   as an holomorphic function satisfying the following
   functional  equation:
 $$\xi(E,\ 2-s) = \pm \xi(E,\ s),$$
 $$\xi(E,\ s) = (2^5pq)^{s/2}(2\pi)^{-s}
\Gamma(s)L(E/{{\bf Q},\ s}).$$ \qquad Our main results are
 to determine Selmer groups, Shafarevich-Tate group and
 Mordell-Weil group(For definitions and notations, see [5]).
For any abelian group $G$, we let $G[n]$ be the $n-$torsion
part of $G$. Take the following elliptic curve $ E^\prime $
and isogeny $\varphi$ of degree 2.
$$\hskip 4cm E^\prime :\; \; y^2 = x^3-\sigma 2(p+q)x^2 + 4x.
 \qquad  \hskip 3.6cm (1.2)  $$
 $$\varphi :\quad E\rightarrow E^\prime,\quad (x,\ y)\mapsto
 (y^2/x^2,\ y(pq-x^2)/x^2).$$
 Then  $\ E[\varphi] =\{O, (0, 0)\}\ $ be the kernel of
 $\varphi.\ $ The dual isogeny of $\varphi $ is
  $$ \widehat{\varphi }:\quad E^\prime \rightarrow E,
  \quad (x,\ y)\mapsto  (y^2/4x^2,\ y(4-x^2)/8x^2).$$

\par \vskip 0cm

     {\bf Theorem 3}.  Let $E=E_+$  be the elliptic curve
     in (1.1). Then on the Selmer groups
 $\; S^{(2)}(E/{{\bf Q}}),\;  $
 $ S^{(\varphi)}(E/{{\bf Q})},\; $  and
  $\; S^{(\widehat{\varphi})}(E/{{\bf Q})};\;  $
    \hskip 0.1 cm
  Shafarevich-Tate group \hskip 0.1 cm
   $\Sh(E/{{\bf Q}} )[2],\; $   Mordell-Weil group
   $\ E({{\bf Q}}),\; $
and rank$E({{\bf Q}})\ $ of $\ E,\ $
we have the following results .  \par \vskip 0.12cm
(1)  \hskip 0.5 cm $S^{(\varphi)}(E/{{\bf Q})}\cong
\left \{ \begin{array}{ll}
\{0\}, &  \mbox{\quad if $p\equiv1, 3, 5\
  ({{\rm mod}} \ 8),$} \\ {{\bf Z}}/2{{\bf Z}},
  &  \mbox{\quad if $p\equiv \; 7\qquad ({{\rm mod}}\ 8).$}
     \end{array} \right. $\par   \vskip 0.12cm

\hskip0.83cm $S^{\widehat{(\varphi)}}(E^\prime /{{\bf Q})}
\cong \left \{
     \begin{array}{ll}
     ({{\bf Z}}/2{{\bf Z}})^2, &  \mbox{\quad
     if $p\equiv\;  5\qquad ({{\rm mod}} \ 8),$} \\
     ({{\bf Z}}/2{{\bf Z}})^3, &  \mbox{\quad
     if $p\equiv 1, 3, 7\  ({{\rm mod}}\ 8).$}
     \end{array} \right. $ \par   \vskip 0.12cm

 (2)  \hskip 0.5cm $S^{(2)}(E/{{\bf Q})}\cong
     \left \{
     \begin{array}{ll}
     ({{\bf Z}}/2{{\bf Z}})^2, &  \mbox{\quad
     if $p\equiv\; 5\quad ({{\rm mod}} \ 8),$} \\
     ({{\bf Z}}/2{{\bf Z}})^3, &  \mbox{\quad
     if $p\equiv 1, 3\  ({{\rm mod}}\
     8),$}\\
     ({{\bf Z}}/2{{\bf Z}})^4, &  \mbox{\quad
     if $p\equiv\; 7\quad  ({{\rm mod}}\ 8).$}
     \end{array} \right. $ \par \vskip 0.12cm

 (3) \hskip 0.26cm
        ${{\rm rank}}E({{\bf Q}}) +
        {\rm dim}_{{{\bf F}}_2}$ $\Sh(E/{{\bf Q}} )[2]
        = \left \{
     \begin{array}{ll}
     0, &  \mbox{\quad if $p\equiv \; 5\quad
     ({{\rm mod}} \ 8),$} \\
     1, &  \mbox{\quad if $p\equiv 1, 3\
     ({{\rm mod}}\ 8),$}\\
     2, &  \mbox{\quad if $p\equiv \; 7\quad
     ({{\rm mod}}\ 8).$}
     \end{array} \right. $  \par

  (4) \hskip 0.22 cm If  $p\equiv 5\
     ({{\rm mod}} \ 8)$ , then \par
   \hskip 5cm    $r={{\rm rank}}E({{\bf Q}}) = 0,$ \par
   \hskip 5cm    $\Sh(E/{{\bf Q}} )[2] =\{0\},$   \par
   \hskip 5cm    $E({{\bf Q}}) \cong
   {{\bf Z}}/2 {{\bf Z}} \oplus
   {{\bf Z}}/2{{\bf Z}}.$
        \par \vskip 0.3cm

   {\bf Theorem 4}.  Let $E=E_-$  be the elliptic curve in
    (1.1). Then on the Selmer
groups, Shafarevich-Tate group, and Mordell-Weil group of $E$,
we have the following results      \par   \vskip 0.16cm
      (1)  \hskip 0.5 cm $S^{(\varphi)}(E/{{\bf Q})}\cong
      \left \{
     \begin{array}{ll}
     \{0\}, &  \mbox{\quad if $p\equiv 3, 5\
     ({{\rm mod}} \ 8),$} \\
     {{\bf Z}}/2{{\bf Z}}, &  \mbox{\quad
     if $p\equiv 1, 7\ ({{\rm mod}}\ 8).$}
     \end{array} \right. $\par \vskip 0.12cm
    \par

\hskip1.2cm $S^{\widehat{(\varphi)}}(E^\prime /{{\bf Q})}
\cong
     ({{\bf Z}}/2{{\bf Z}})^2. $ \par  \vskip 0.12cm

     (2)  \hskip 0.5cm $S^{(2)}(E/{{\bf Q})}\cong
     \left \{
     \begin{array}{ll}
     ({{\bf Z}}/2{{\bf Z}})^2, &  \mbox{\quad
     if $p\equiv 3, 5\ ({{\rm mod}} \ 8),$} \\
     ({{\bf Z}}/2{{\bf Z}})^3, &  \mbox{\quad
     if $p\equiv 1, 7\  ({{\rm mod}}\
     8).$}
     \end{array} \right. $ \par \vskip 0.12cm
     (3)  \hskip 0.26cm
        ${{\rm rank}}E({{\bf Q}}) +
         {\rm dim}_{{{\bf F}}_2}$ $\Sh(E/{{\bf Q}} )[2]
          = \left \{
     \begin{array}{ll}
     0, &  \mbox{\quad if $p\equiv 3, 5\
     ({{\rm mod}} \ 8),$} \\
     1, &  \mbox{\quad if $p\equiv 1, 7\
     ({{\rm mod}}\ 8).$}
     \end{array} \right. $  \par \vskip 0.12cm
     (4)  \hskip 0.22 cm If  $p\equiv 3, 5\
     ({{\rm mod}} \ 8)$ , then \par
   \hskip 5cm    $r={{\rm rank}}E({{\bf Q}}) = 0,$ \par
   \hskip 5cm    $\Sh(E/{{\bf Q}} )[2] =\{ 0 \} ,$   \par
   \hskip 5cm    $E({{\bf Q}}) \cong
   {{\bf Z}}/2 {{\bf Z}} \oplus
   {{\bf Z}}/2{{\bf Z}}.$
                  \vskip 0.38cm

 {\bf Theorem 5}. If $\ p\equiv 3\ ({{\rm mod}}\ 8)\; $
 and
 $\; q=a^2+b^2,\; (a+\varepsilon)^2+(b+\delta)^2=c^2\ $ for
 some rational integer $\ a,\ b,\ c\ $ and
 $\ \varepsilon , \delta=\pm 1,\; $
 then  the elliptic curve $E=E_+$
in (1.1) has rank
$$r={{\rm rank}}E({Q}) = 1 .$$

      {\bf Example 1}.  Let $(p,\ q) = (3,\ 5),\ $
consider the elliptic curve  $\ E :\;  y^2 = x (x+3)(x+5).\; $
       By Theorem 5 we have
 $\ r={{\rm rank}}E({{\bf Q}}) =1.\; $
 Thus $\ E({{\bf Q}}) \cong \{ O, \;  (0,\ 0), \;
        (-3,\ 0), \; (-5,\ 0) \} \oplus {{\bf Z}}.\;  $
  We could also find a point in $ E({{\bf Q}}) $
  with infinite order: $\ P=(-4 ,\ 2) $.
   \par \vskip 0.12cm

      {\bf Example 2}.  Let $(p,\ q) = (11,\ 13),\ $
  consider the elliptic curve $E :\;  y^2 = x (x+11)(x+13).\; $
   By Theorem 5 we have
  $r={{\rm rank}}E({{\bf Q}}) =1.\; $
   Thus $E({{\bf Q}}) \cong \{ O, \; (0,\ 0), \;
   (-11,\ 0), \;  (-13, \ 0) \} \oplus {{\bf Z}}.\; $
   By handcount, we find a point of
infinite order in $E({{\bf Q}}): \; $
$$ P =(243391201/1587600, \;
4094288981999/2000376000). $$
 \par \vskip 0.66cm

\begin{center}
{\Large II.\ Proofs of Theorem 1 and 2}
\end{center}

First calculate quantities associated with the elliptic curve
$E$ in (1.1) (for notations see [5]):
 $a_1=a_3=a_6=0,\ a_2 =\sigma (p+q),\ a_4=pq,$
  $b_2=\sigma 4(p+q),\ b_4=2pq,\ b_6=0,\ b_8=-p^2q^2, $
  $\ c_4=16(p^2+2q),\ c_6=\sigma32(p+q)(pq-8),\ $
$\Delta=\Delta(E)=64p^2q^2,\; j=j(E)=64(p^2+2q)^3/p^2q^2.$
\par \vskip 0.28cm

{\bf Lemma 2.1}
(1) The equation (1.1) for $E=E_\sigma $ is a global minimal
Weierstrass-equation.

(2) The reduction of $E$ is additive at 2,  multiplicative
at $p$ and $ q$, and good at other primes.

(3)  $\ E_+ \ $ has a split multiplicative reduction at
$p$ if and only if it does at $q$, which is just equivalent to
 $\ (2/p)=1,\; i.e.,\;  p\equiv 1,\ 7 \ (mod\ 8) $.

(4)  $E_- \ $    has
 a split multiplicative reduction at $p$ if and only if
$ (-2/p)=1,\ i.e.,\  p\equiv 1,\ 3 \ (mod\ 8) ;\; $
 $E_-  $ has such reduction at $q$  if and only if  $
(2/q)=1,\ i.e.,\ $  $  q\equiv 1,\ 7 \ (mod\ 8) $.
\par \skip 0.4cm

{\bf Proof.} Let $v_\ell$ denote the normalized $\ell-$adic
exponential valuation of
{${\bf Q}$} ($\ell$ is a prime number).
(1) Since $v_\ell (\Delta (E))= v_\ell
(64p^2q^2)<12$ for every $\ell$, so
equation (1.1) is minimal at every prime $\ell $, i.e.,
 it is global minimal.
 (2) Note that
$\ v_\ell (\Delta)=0\ $ when $\ell \neq 2,\ p,\ q$.
 And $v_2(c_4)>0, \ v_p(c_4)=v_q(c_4)=0. $
 The assertion follows from Prop.5.1 of [5, p180].
(3) The reduction of $E$ at $p$ is $y^2=x^2(x+2)$,
which has unique singular point at
 $(0,\ 0)$ with ``tangent" lines $\ y^2=2x^2, \ $ or
  $\ y=\pm \sqrt{2}x\ $. Thus the reduction is split if and
 only if 2 is a square residue modulo $p,\ $ i.e.,
  $\ (2/p)=1.\ $ The same thing happens for $q$.
 (4) The proof is similar to that of (3).
 \par \vskip 0.3cm

{\bf Lemma 2.2.} The elliptic curve $E$ in (1.1) is
supersingular at prime $\ell \ (\neq 2,
p, q)$ if and only if
$$\sum_{k=0}^{m}(C_m^k)^2p^kq^{m-k} \equiv 0\; (mod \ \ell).
\qquad ({{\rm here}} \quad  m=(\ell-1)/2 \ )$$
{\bf Proof.} Let $y^2=f(x)=x(x+\sigma p)(x+\sigma q)$.
The coefficient of $\ x^{\ell-1}\ $ in  $\ f(x)^m\ $ is
$$\sum_{i+j=m}C_m^{m-i}C_m^{m-j}(\sigma p)^{i}(\sigma q)^{j}=
 \sigma^m\sum_{k=0}^m(C_m^k)^2p^kq^{m-k}.$$
Now the lemma follows from Theorem 4.1 of [5, p140].
\par \vskip 0.3cm

{\bf Lemma 2.3}$^{[9]}$. Let $E$ be any elliptic curve
  over a number field $K$, $\wp $
a prime ideal of $K$ over a prime number
 $\ell $, $k_\wp $ the residue class field of $K$
modulo $\wp$, $\ e = e(\wp |\ell) $ and $ f = f(\wp |\ell) $
 the ramification index and
residue class degree respectively of
 $\wp$ over $\ell$,  $\ \tilde{E}$
the reduced curve over $k_\wp$ of $E$ modulo $\wp$.
Put $t_\ell =0 $ (if $\ell-1>e$)\ or \ max $\{ r \ \ :\
  (\ell-1)\ell ^{r-1}\leq e,\ 0<r\in {\bf Z}\} $
(otherwise).

(1) If $E$ has good reduction  at $\wp$, then
 $$\#E(K)_{tors}\ \mid \  \#\tilde{E}(k_\wp)\  \ell^{2t_\ell}
\ \leq \ (1+\ell^{f}+2\ell ^{f/2})\ \ell ^{2t_\ell}.$$
\hskip 0.9cm(2) If $E$ has additive reduction  at $\wp$,
 then
$$\#E(K)_{tors}\ | \ 12 \ \ell ^{2(t_\ell+1)}.$$

{\bf Proof of Theorem 1}. By Lemma 2.1 we know $E$ has a
 good reduction $\tilde{E}$ at 3.
Now $m=(3-1)/2=1,\ $ so
 $\ \sum_{k=0}^{m}(C_m^k)^2p^kq^{m-k}=p+q\equiv (-1)+(1)
 \equiv 0\ ({{{\rm mod}}}\ 3) $
since $p, q$ are twin prime numbers,
 which means, by Lemma 2.2, the reduced elliptic
curve   $\tilde{E}$ over ${{\bf F}}_3$
is supersingular. In particular, we know
\#$\tilde{E}({{\bf F}_3)} = 3+1 = 4\ $
 (See [5], p145, Exer. 5.10(b)).

  Now, as an elliptic curve over $K,\ $ $E$  obviously
 has a good reduction at $\wp$ since $\wp | 3$.
 If $f=f(\wp |3)=1$, then
  $k_\wp = {{\cal O}}_K/\wp =
  {{\bf Z}}/(3) = {{\bf F}}_3,\
 $ so the reduced curves of $E$ modulo $\wp $ and
  modulo $3$ is the same curve
 ${\tilde{E}}$, which is supersingular. Thus
  $\ \#\tilde{E}(k_\wp) = \#\tilde{E}({{\bf F}}_3) = 4.\ $
 By Lemma 2.3 we know $\ \#E(K)_{tors} \ |\
  \#\tilde{E}(k_\wp)\ 3^{2t_3} = 4\cdot 9^{t_3},\
 $ where $t_3 $ is as in Lemma 2.3. Also since
 $E(K)[2]=E[2]\subset E(K)_{tors},\ $ so
  $E(K)_{tors}/E(K)_3 \cong E(K)[2] \cong
  {{\bf Z}}/2{{\bf Z}}\oplus
  {{\bf Z}}/2{{\bf Z}}.\ $
  (In fact we have
  ord$_3(\#E(K)_3) \leq 2\ (1+{{\rm log}}_3(e/2))\ $
   by Lemma 2.3.)   This proves (3).

Now if $f=e=1,\ $ then $3-1>e$, by Lemma 2.3 we know $t_3=0$,
  which proves (1).

  If $E$ has an additive reduction at a prime $\wp$ of $K$,
  then by  Prop.5.4(b) of [5, p181],
   we know  $E$ has  additive reduction at  prime  number
   $\ell $ under $\wp$ , and by Lemma 2.1 we know $\ell =2$,
   so  $\wp | 2.\ $ Thus by Lemma 2.3(2) we know
   $\#E(K)_{tors}$ divides $12\cdot 2^{2(t_2+1)}.\ $
   So  $\#E(K)_{tors}$ only has prime-factors 2 and 3, and
   the order of 3 in it is 0 or 1.  This proves (4).

  Combining (3) and (4) gives (2).  This proves Theorem 1.
 \par \vskip 0.3cm

   {\bf Proof of Theorem 2}.
   By Lemma 2.1 and the definition of $f_\ell$ we know that
$f_\ell=2+\delta_2;\ 1 ;\ $ or $\ 0\ $ respectively if
 $\ell =2;\ \ell=p$ or $q;\ $ or otherwise.
 We now use Tate's method in [7] to determine $\delta_2$ and
the type of Kodaira symbol. Make the isomorphic transformation
 $ \rho:\; E\rightarrow E^\prime, \quad x\mapsto
 x^\prime+\sigma,\; y\mapsto y,\ $ $(\sigma =\pm 1)$
  then the elliptic curve
  $E^\prime $ has the Weierstrass equation
 $$ E^\prime :\; {y^\prime }^2 =
 (x^\prime +\sigma)\ (x^\prime +\sigma +\sigma p) \
 (x^\prime +\sigma +\sigma q).$$
 It is easy to see that the constants
 associated to $E^\prime$  (defined as in [5, p46]) satisfy
$ \Delta^\prime=\Delta (E^\prime)=64p^2q^2 $
  $ \equiv a_3^\prime \equiv a_4^\prime \equiv a_6^\prime
   \equiv b_2^\prime
\equiv 0 $ $ \ ({{{\rm mod}}}\ 2) ,\; $
$ a_6^\prime \equiv b_6^\prime \equiv b_8^\prime \equiv 0 (
{{{\rm mod}}}\ 4);\; $
$b^\prime_8 \equiv 4 \not\equiv 0 $ (mod 8).
 According to Tate's criterion in [6-7],
 we know that the type of the Kodaira symbol of $E$ is III,
 and the local  exponent $ f_\ell $  of conductor $ N_E $ is
 $f_2 = v_2(\Delta (E))-1=v_2(64p^2q^2)-1=5,\ $  so
 $\delta_2 = 5-2=3 $.
 By Ogg's formula [5, p361] we obtain that
 the number of irreducible components on the special
fiber ${{\cal C}}_\ell$ is
 $m_\ell = v_2(\Delta(E))+1-f_2=6+1-5=2.$

 Then consider the case of prime $p$. Since
 $ \Delta (E) \equiv a_3\equiv a_4\equiv a_6 \equiv 0 \
 ({{{\rm mod}}}\ p),\ b_2 \not \equiv 0
\; ({{{\rm mod}}}\ p),\ $ $v=v_p(\Delta (E))=2,\ $
 so by Tate's criteria in [7] we know that the Kodaira symbol
  of $E$ at $p$ is $I_v = I_2$.
  If $E$ has a non-split multiplicative reduction at $p$
  ,  then the Tamagawa constant $c_p=2$
  since $v_p(\Delta) =2$ is even. And if $E$ has a split
   multiplicative reduction at $p$, then $c_p=v_p(\Delta) =2$
 by the Kodaira-Neron theorem [5, p183, Thm. 6.1].
 Thus $E({{\bf Q}}_p)/E_0({{\bf Q}}_p)
 \cong {{\bf Z}}/2{{\bf Z}}$ as desired.
 The case for $q$ goes similarly.
 \par \vskip 0.12cm

 {\bf Diamond's Theorem}$^{[8]}$. An elliptic
 curve $E$ defined over rationals ${{\bf Q}}$
 is modular if it has good or multiplicative reduction at
 3 and 5.
 \par \vskip 0.12cm

 {\bf Proof of Corollary 2.1}.  By Lemma 2.1 we know
$E$ has good or multiplicative reduction at
 3 and 5, so $E$ is modular via Diamond's theorem.
 Thus Hasse-Weil conjecture is true for $E$.
\par \vskip 0.76cm

\begin{center}
{\Large III. Proof of Theorem 3 and 4}
\end{center}

 Denote by $M_{{\bf Q}}$  the set of primes
 (or normalized valuations) of the
 rational field ${{\bf Q}}$ containing the
  infinity prime $\infty $.
Let $S=\{\infty ,\ 2,\ p,\ q\}$,  and define
\begin{center}
${{\bf Q}}(S, 2)=\{d\in{{\bf Z}} : \ d  | 4 $
is squarefree $\}=\{\pm 1,\ \pm 2\}\; $
(regarded as subgroup of
  ${{\bf Q}}^*/{{{\bf Q}}^*}^2$);
\end{center}
 \begin{center}
 ${{\bf Q}}(S, 2)^\prime= \{ d\in{{\bf Z}} :\ d
   |  pq \ $ is squarefree $ \} =
  \{ \pm 1,\ \pm p,\ \pm q,\ \pm pq \}\; $
  (as subgroup of ${{\bf Q}}^*/{{{\bf Q}}^*}^2 $).
   \end{center}
For $d_1\in {{\bf Q}}(S, 2),\; d_2=4/d_1,\  $
  define the curve
  $$ C_{(d_1)}: \;  y^2 =d_1 x^4 -\sigma 2(p+q)x^2 + d_2. $$
For $d_1\in {{\bf Q}}(S, 2)^\prime,\; d_2=pq/d_1,\  $
define the curve
 $$C^\prime_{(d_1)}:\; y^2 =d_1x^4+\sigma (p+q)x^2+d_2.$$
  Then according to Cremona's algorithm in [11, p63-65],
   we have

  {\bf Lemma 3.1}.
  $S^{(\varphi )}(E/{{\bf Q}})\ \cong \
  \{d_1\in{{\bf Q}}(S, 2)\ : \
C_{(d_1)}({{\bf Q}}_v) \neq \emptyset ,\; \forall v\in S\}
\ \subset \ \{\pm 1,\ \pm 2 \};$
\qquad $$S^{(\widehat{\varphi })}(E^\prime /{{\bf Q}})
\ \cong \ \{d_1\in{{\bf Q}}(S, 2)^\prime\ : \
C^\prime_{(d_1)}({{\bf Q}}_v) \neq \emptyset ,\;
 \forall v\in S \}
\ \subset \ \{\pm 1,\ \pm p,\ \pm q,\ \pm pq\}. $$
\qquad{\bf Lemma 3.2} (Hensel's lemma). Suppose that $f\in
R[x_1,\ \cdots ,\ x_n] $ is a polynomial over a
  complete ring $R$ with valuation $v,\ $  If
  there exists  $\ a\in R^n\ $ such that
 $v(f)>2v(f_{k} (a)) \; $
 (for any $k$ with $ 1\leq k \leq n),\; $
   then $f$ has a root in $R^n.$
 (See [5, p322, Ex.10.12], Here $f_{k}=f_{x_k}={\partial
f}/{\partial x_k}$.)

 {\bf Lemma 3.3}$^{[12]}$. A unit $\alpha \in {{\bf Q}}_2$
 is a square in ${{\bf Q}}_2$ if and only if
 $v_2(\alpha -1)\geq 3 $.

 {\bf Lemma 3.4}. Let  ${{\bf Z}}_\ell $
 be the ring of $\ell-$adic integers,
$\ell $ a prime number. Consider    $y^2=f(x)$  with
  $f\in {{\bf Z}}_\ell [x],\
 {\rm deg} f =4,\  $  discriminate
 ${\rm disc} (f)\neq 0.\; $
 Denote by $\bar{f} $  the image of $f$ under the map
 ${{\bf Z}}_\ell [x]\rightarrow {{\bf F}}_\ell [x]$
 deduced from the natural map
 ${{\bf Z}}_\ell \rightarrow {{\bf F}}_\ell . $
 If $\bar{f}$ is not a constant
modulo a square (i.e.,
$\ \bar{f}=\bar{g}^2\bar{h} \ $
 with   deg$\bar{h}\geq 1,\ $  deg$\bar{g}\geq 0,\ $
 $ \bar{h}$ is squarefree,
  $\ \bar{g}, \bar{h}\in {{\bf F}}_\ell [x]),\ $
  then $y^2=f(x)$ has a solution in
 ${{\bf Q}}_\ell$. ( See [10, p140, Lemma 14.])
   \par \vskip 0.2cm

 {\bf Proposition 3.5}. Let $C=C_{(2)}:\
 y^2=f(x)=2x^4-\sigma 2(p+q)x^2+2 $ be one of the
  $  C_{{(d_1)}}$.

 (1) If
 $p\equiv 7\ ({{\rm mod}}\ 8)$,
 then $C_{(2)}({{\bf Q}}_2)\neq \emptyset .$

 (2) $C_{(2)}({{\bf Q}}_p)\neq \emptyset $ if and only if
 $(2/p)=1, \; i.e.,\; p\equiv 1,\ 7\ ({{\rm mod}}\ 8).$

 (3) $C_{(2)}({{\bf Q}}_q)\neq \emptyset $ if and only if
 $(2/q)=1, \; i.e.,\; q\equiv 1,\ 7\ ({{\rm mod}}\ 8)$
 (that is $ p\equiv 5,\ 7\ ({{\rm mod}}\ 8)$)
\par \vskip 0.4cm

{\bf Proof}. We prove for $\sigma=1$.
(1) Let $g(x, y)=f(x) -y^2$, then
 $g(1,\ -2)=-4(p+1),\ g_y(1,\ -2)=4$.
If $p\equiv 7 \ ({{\rm mod}}\ 8)$,
 then
$v_2(g(1,\ -2))  > 2 v_2(g_y(1, -2))$.
 Hensel's  lemma  tells us
$C({{\bf Q}}_2)\neq \emptyset .$

(2) If $C({{\bf Q}}_p)\neq \emptyset $, we assert that
$C({{\bf Z}}_p)\neq \emptyset. $  In fact, if
 $(x ,\ y )\in C({{\bf Q}}_p),\ $
  then by the equation obviously $(1/x,\ y/x^2)$
 is also in $C({{\bf Q}}_p)$. So we have either
   $v_p(x)\geq 0$ (and then $ v_p(y)\geq 0 $), or else
$v_p(x)<0$ (and then $ v_p(1/x)>0$ and $v_p(y/x^2)\geq 0$),
our assertion follows.
Now assume  $(x ,\ y )\in C({{\bf Z}}_p),\ $
so $\ y^2=f(x)= 2(x^2-1)^2 - 4px^2.\ $
Thus $y^2\equiv 2(x^2-1)^2\ ({{\rm mod}}\ p),\ $
which means $(2/p)=1$ (Note
that $ y\not \equiv 0 \ ({{\rm mod}}\ p),\ $
 otherwise we'd have $x^2-1\equiv 0 \
({{\rm mod}}\ p),\ $
$ x \not \equiv 0\ ({{\rm mod}}\ p),\  $
so $v_p(y^2)=v_p(4px^2)=1,\ $
  a contradiction). Conversely, if $(2/p)=1, $  then
$2\equiv a^2\ ({{\rm mod}}\ p), $ $a\in {{\bf Z}}$.
Then obviously $v_p(g(p, a))>0,\
v_p(g_y(p,a))=v_p(-2a)=0$. So
$C({{\bf Q}}_p)\neq \emptyset $ by Hensel's lemma.

(3) The proof is similar to that of (2).
\par \vskip 0.26cm

{\bf Proposition 3.6}. For the curve
$\; C^\prime=C^\prime_{(-1)}:\;
 y^2=f(x)=-x^4+(p+q)x^2-pq $, we have

 \quad (1) $C^\prime({{\bf Q}}_p)\neq \emptyset;\quad $
 $\ C^\prime({{\bf Q}}_q)\neq \emptyset ; $

\quad (2) $C^\prime ({{\bf Q}}_2)\neq \emptyset\; $
if and only if
 $\; p \equiv 1,\ 3, \ 7\ ({{\rm mod}}\ 8).$
\par \vskip 0.2cm

{\bf Proof}. (1) Note $f = -(x^2-p)(x^2-q),\ $ so
$f\equiv x^2(-x^2+2)\ ({{\rm mod}} \ p).\; $
Since $-x^2+2 $ is not a square in
${{\bf F}}_p[x]$,
 Lemma 3.4 tells us
$C^\prime({{\bf Q}}_p)\neq \emptyset. $
The case for $q$ is Similar.

(2) Consider the sufficiency first. (i) If
$p \equiv 1 ({{\rm mod}}\ 8),\ $ then
$p=\alpha ^2$ for some $\alpha \in {{\bf Q}}_2$
(by Lemma 3.3),
then it is easy to see $(\alpha ,\ 0)$
 is a point of $C^\prime $.
(ii) If $p \equiv  3, \ 7\ ({{\rm mod}}\ 8),\ $ then
let $g(x, y)=y^2-f,\ $ so
$g(0,\ p) =pq+p^2=p(p+q)=2p(p+1),\; g_y(0,\ p)=2p.\ $
So $v_2(g(0, p))=v_2(2p(p+1))\geq 1+2 =3 >
2v_2(g_y(0,\ p))$. Thus
$C^\prime ({{\bf Q}}_2)\neq \emptyset $
 by Hensel's lemma.

For necessity, assume
$(a, \ b)\in C^\prime ({{\bf Q}}_2) $, then
obviously $b^2=1-(a^2-p-1)^2$.
 If $b$ is not a 2-adic integer, then $a$
too; assume $a=2^{-n} a_0,\ b=2^{-2n} b_0
,\ v_2(a_0)=v_2(b_0)=0  $
(since $v_2(b^2)=2v_2(a^2-p-1)=2v_2(a^2));\ $
so $2^{-4n} b_0^2=1-(2^{-2n} a_0^2-p-1)^2,$
$\ b_0^2+a_0^4= 2^{2n}(2^{2n}+2(p+1)a_0^2-2^{2n}(p+1)^2)
\equiv 0 \ ({{\rm mod}}\ 4)$, contradict to
 $b_0^2+a_0^4\equiv 1+1 \ ({{\rm mod}}\ 8)$ by
Lemma 3.3. Thus $b\in {{\bf Z}}_2$.
 Then, since $b^2=1-(a^2-p-1)^2=1-c^2 \ $,
 there are only two possibilities:
 $(b^2,\ c^2)\equiv (1,\ 0)  $ or
  $(0,\ 1)\ ({{\rm mod}}\ 8)$.
(i) $c^2\equiv 0\ ({{\rm mod}}\ 8)$ means
$p\equiv a^2-1\equiv 0-1 \ ({{\rm mod}}\
4).\ $ (ii) $c^2\equiv 1\ ({{\rm mod}}\ 8)$ means
 $ a\equiv 1 \ ({{\rm mod}}\ 2). \ $  By  Lemma 3.3
 we know $b^2=1-c^2\equiv 0 \ ({{\rm mod}}\ 16)$.
  Rewrite $b^2=1-c^2$ as $b^2=1-p^2-(a^2-1)^2+2p(a^2-1)$,
   then by Lemma 3.3 we know $0\equiv b^2\equiv
1-p^2+0+0 \ ({{\rm mod}}\ 16),\ $ that is
$p\equiv \pm 1 \ ({{\rm mod}}\ 8).\ $
This proves Proposition 3.6.
\par \vskip 0.3cm

{\bf Proof of Theorem 3}.

(1) Consider $S^{(\varphi  )}(E/{{\bf Q}})$.
Look at  Lemma 3.1 and definitions before it,  if
$d_1\in {{\bf Q}}(S,\ 2) $ and $d_1<0$, then $d_2=4/d_1<0$
; then obviously  $C_{(d_1)} $  has no
solution in ${{\bf R}} = {{\bf Q}}_\infty$;
 so we must have $d_1\geq 0,\ i.e.,\ $
 $S^{(\varphi  )}(E/{{\bf Q}})\cong \{1 \}$ or
$\{1,\ 2 \}. $ Also obviously
 $(0,\ \sqrt{2})\in C_{(2)}({{\bf R}}),\ $
so $2\in S^{(\varphi  )}(E/{{\bf Q}})\Leftrightarrow
 C_{(2)}({{\bf Q}}_v)\neq \emptyset\ $
 (for $v=2, p, q $) $\Leftrightarrow
 p\equiv 7   \ ({{\rm mod}}\ 8)\  $
 (The last step is via Proposition 3.5).
Thus $S^{(\varphi  )}(E/{{\bf Q}})\cong (0)$
 (if $p\equiv 1,3,5  \  ({{\rm mod}}\ 8)) $ or
 ${{\bf Z}}/2{{\bf Z}}$ (if
$p\equiv 7 \ ({{\rm mod}}\ 8) $).

As for $S^{(\widehat{\varphi  })}(E^\prime/{{\bf Q}}), $
obviously  $C^\prime_{(-p)}$ and $C^\prime_{(-q)}$ both
have rational point $(1,  0)$, so
$-p,-q\in S^{(\widehat{\varphi  })}(E^\prime /{{\bf Q}})$,
 $\{1,\ -p,\ -q,\ pq\}\subset
 S^{(\widehat{\varphi  })}(E^\prime /{{\bf Q}}).$
Note that $(\sqrt p,\ 0)$
is a real point of $C^\prime_{(-1)}$,
so, by Proposition 3.6,
$\ -1\in S^{(\widehat{\varphi  })}(E^\prime /{{\bf Q}})
\Leftrightarrow p\equiv 1, 3, 7   \ ({{\rm mod}}\ 8).\ $
This gives $S^{(\widehat{\varphi  })}(E^\prime /{{\bf Q}})$
as desired.

(3) Since $E^\prime ({{\bf Q}})[\widehat{\varphi}]
 =\{O,\ (0,0)\}=\varphi (E({{\bf Q}})[2]),\ $
 so by the exact sequence (see [5, p301]) :
 $$ 0\rightarrow \frac{E^\prime ({{\bf Q}})
 [\widehat{\varphi}]} {\varphi (E({{\bf Q}})[2])}
 \rightarrow
 \frac{E^\prime({{\bf Q}})} {\varphi(E({{\bf Q}}))}
\stackrel{\widehat{\varphi}}{\longrightarrow}
 \frac{E({{\bf Q}})}{2E({{\bf Q}})}\rightarrow
 \frac{E({{\bf Q}})}
{\widehat{\varphi}(E^\prime({{\bf Q}}))} \rightarrow 0 ,$$
 we have the isomorphism of
 vector spaces over the finite field ${{{\bf F}}}_2 :$
$${E({{\bf Q}})}/{2E({{\bf Q}})}\cong
 \frac{E^\prime({{\bf Q}})} {\varphi(E({{\bf Q}}))}
 \oplus
\frac{E({{\bf Q}})} {\widehat{\varphi}
(E^\prime({{\bf Q}}))} . $$
On the other hand, denote
$r={\rm rank}\ E({{\bf Q}}), \ $  we have
$$\frac{E({{\bf Q}})}{2E({{\bf Q}})}\cong
 ({{\bf Z}}/2{{\bf Z}})^{r}\oplus
E({{\bf Q}})[2]\cong ({{\bf Z}}/2{{\bf Z}})^{r+2},$$
 So we have the following dimension formula
(here ${{\rm dim}} =
{{\rm dim}}_{\stackrel{}{{{\bf F}}_2}} $
) : $$r={{\rm dim}} \frac{E^\prime({{\bf Q}})}
{\varphi(E({{\bf Q}}))} +
{{\rm dim}} \frac{E({{\bf Q}})}
{\widehat{\varphi}(E^\prime({{\bf Q}}))} -2 . $$
From the following exact sequences (and their duals)
(see [5, p298, 314]):
  \begin{center}
$0\rightarrow
{E^\prime({{\bf Q}})} /{\varphi(E({{\bf Q}}))}
\rightarrow S^{({\varphi  })}(E/{{\bf Q}})
\rightarrow $  $\Sh(E/{{\bf Q}} )[\varphi]\rightarrow 0 ;\; $
     \end{center}
    \begin{center}
$0\rightarrow $   $\Sh(E/{{\bf Q}} )[\varphi]\rightarrow $
  $\Sh(E/{{\bf Q}} )[2]
 \stackrel{\widehat{\varphi}}{\longrightarrow}  $
  $\Sh(E^\prime/{{\bf Q}} )[\widehat{\varphi}]\rightarrow 0, $
        \end{center}
     we have
    \begin{center}
$ S^{({\varphi  })}(E/{{\bf Q}})\cong
{E^\prime({{\bf Q}})}/ {\varphi(E({{\bf Q}}))}\oplus
 $  $\Sh(E/{{\bf Q}} )[\varphi] ; $ \end{center}
   \begin{center}
$\; S^{(\widehat{\varphi  })}(E^\prime/{{\bf Q}})\cong
{E({{\bf Q}})}/ {\widehat{\varphi}
(E^\prime({{\bf Q}}))}\oplus
$  $\Sh(E^\prime/{{\bf Q}} )[\widehat{\varphi}];$
                 \end{center}
    \begin{center}
 $\Sh(E/{{\bf Q}} )[2]\cong $
 $\Sh(E/{{\bf Q}} )[\varphi]\ \oplus\  $
 $\Sh(E^\prime/{{\bf Q}} )[\widehat{\varphi}]. $
                       \end{center}
Take together we have
    \begin{center}
$r= {{\rm dim}} S^{({\varphi  })}(E/{{\bf Q}}) -
{{\rm dim}}$   $\Sh(E/{{\bf Q}} )[\varphi] +
{{\rm dim}}S^{(\widehat{\varphi  })}(E^\prime/{{\bf Q}})
 - {{\rm dim}}$
  $\Sh(E^\prime/{{\bf Q}} )[\widehat{\varphi}] -2 $
             \end{center}
    \begin{center}
 $={{\rm dim}} S^{({\varphi  })}(E/{{\bf Q}}) +
{{\rm dim}}S^{(\widehat{\varphi  })}(E^\prime/
{{\bf Q}}) - {{\rm dim}}$
 $\Sh(E/{{\bf Q}} )[2] -2 .$
\end{center}
 Then by the result in (1) we obtain
 $r+ {{\rm dim}}$  $\Sh(E/{{\bf Q}} )[2]=0,\ 1,\ 2 $
 in cases as in Theorem 3(3).

(2) From the exact sequence ([5, p305]) :
    \begin{center}
$0\rightarrow
{E({{\bf Q}})}/ {2E({{\bf Q}})}
\rightarrow S^{({2})}(E/{{\bf Q}})
\rightarrow $  $\Sh(E/{{\bf Q}} )[2]\rightarrow 0 , $
        \end{center}  we have
$ S^{({2})}(E/{{\bf Q}})\cong
{E({{\bf Q}})}/ {2E({{\bf Q}})}\oplus
 $  $\Sh(E/{{\bf Q}} )[2] , $
 thus obtain the dimension formula
${{\rm dim}}S^{({2})}(E/{{\bf Q}})=2;\ 3;\  4\ $
  if $p\equiv 5;\ 1 $ or $3;\ 7 \ ({{\rm mod}}\ 8 ) $
  respectively, from which we  obtain easily the
  isomorphic type of $ S^{({2})}(E/{{\bf Q}}) $
  as Theorem  3(2) asserted.

  (4) By (3) we now have
  $r = {{\rm dim}}$  $\Sh(E/{{\bf Q}} )[2] = 0,\ $
  so  $\Sh(E/{{\bf Q}} )[2] = \{0\}, \ $ and
  $E({{\bf Q}} ) = E({{\bf Q}} )_{tors}
  \cong {{\bf Z}}/2{{\bf Z}}
  \oplus{{\bf Z}}/2{{\bf Z}}$ by theorem 1.
  This proves Theorem 3.

\vskip 0.5cm
Now consider the proof of Theorem 4, then
we have $\sigma =-1, \ E=E_\sigma = E_- $.
 We  use the same formulae
 as above (but $\sigma =-1 $) to  define
set $S$,
 ${{\bf Q}}(S,\ 2),\; $
 ${{\bf Q}}(S,\ 2)^\prime , \; $
  $ C_{(d_1)} ,\ $ and
$C^\prime_{(d_1)} . $
  Then it is easy to show that Lemma 3.1
  and Proposition 3.5 are still valid (for $\sigma =-1$).
  In addition we
  could  prove in a similar way the following propositions.

  {\bf Proposition 3.7}.  Consider the curve $\ C_{(-2)}:\
 y^2=-2x^4+2(p+q)x^2-2 $ .

 \quad (1) If
 $p\equiv 1\ ({{\rm mod}}\ 8)$,
 then $C_{(-2)}({{\bf Q}}_2)\neq \emptyset .$

 \quad (2) $C_{(-2)}({{\bf Q}}_p)\neq \emptyset $
 if and only if
 $(-2/p)=1, \ i.e.,\ p\equiv 1, 3\ ({{\rm mod}}\ 8).$

 \quad (3) $C_{(-2)}({{\bf Q}}_q)\neq \emptyset $
 if and only if
 $(-2/q)=1, \ i.e.,\ q\equiv 1, 3, \;  $
  $ p\equiv 1, 7\ ({{\rm mod}}\ 8).$
\par \vskip 0.3cm

{\bf Proposition 3.8}.  Consider the curve $C_{(-1)}:\
 y^2=-x^4+2(p+q)x^2-4 $ .

 \quad (1) $C_{(-1)}({{\bf Q}}_p)\neq \emptyset $
 if and only if
 $(-1/p)=1, \ i.e.,\ p\equiv 1\  ({{\rm mod}}\ 4).$

 \quad (2) $C_{(-1)}({{\bf Q}}_q)\neq \emptyset $
 if and only if
 $(-1/q)=1, \ i.e.,\ q\equiv 1\ ({{\rm mod}}\ 4).$
\par \vskip 0.3cm

{\bf Proof of Theorem 4}. By Lemma 3.1 and Proposition 3.5
we obtain
$\ 2\in S^{(\varphi
)}(E/{{\bf Q}})\Leftrightarrow
 p\equiv 7   \ ({{\rm mod}}\ 8).\ $
Lemma 3.1 and Proposition 3.7 give us
 $\ -2\in S^{(\varphi  )}(E/{{\bf Q}})\Leftrightarrow
 p\equiv 1   \ ({{\rm mod}}\ 8).\ $
And Lemma 3.1 and Proposition 3.8 tell us
 $\ -1\not \in S^{(\varphi  )}(E/{{\bf Q}}).\  $
Therefore we have  $\ S^{(\varphi  )}(E/{{\bf Q}})\cong
 \{1\}{{{\bf Q}}^*}^2 ;\;
 \{1, -2\}{{{\bf Q}}^*}^2;\;
  {{\rm or}}\; \{1, 2\}{{{\bf Q}}^*}^2 \; $
 if $\ p\equiv
3, 5;\ 1;\ {{\rm or}}\ 7 \ ({{\rm mod}}\ 8) $
 respectively. That is,
$S^{(\varphi  )}(E/{{\bf Q}})\cong (0)\ $ (if
$\ p\equiv 3, 5 \ ({{\rm mod}}\ 8) )\ $
or $\ {{\bf Z}}/2{{\bf Z}}\ $ (if
 $p\equiv 1, 7 \ ({{\rm mod}}\ 8) $) as desired.

Consider  $S^{(\widehat{\varphi  })}(E^\prime/{{\bf Q}}).$
Since $C^\prime_{(p)}$ and $C^\prime_{(q)}$ both contain
  point $(1, \ 0)$, so
$p, q\in S^{(\widehat{\varphi  })}(E^\prime/{{\bf Q}})$,
 $\{1,\ p,\ q,\ pq\}\subset
 S^{(\widehat{\varphi  })}(E^\prime/{{\bf Q}}).$
On the other hand, if $d_1\in {{\bf Q}}(S,\ 2)^\prime $
 is negative, then $d_2$ is
negative, $C^\prime _{{d_1, d_2}}$ has no real point,
 so $d_1\not\in S^{(\widehat{\varphi
})}(E^\prime/{{\bf Q}}) $.
 Therefore $S^{(\widehat{\varphi  })}(E^\prime/{{\bf Q}})
\cong E({{\bf Q}})/
\widehat{\varphi}(E^\prime ({{\bf Q}})) \cong
\{1, \ p,\ q,\ pq\}{{\bf Q}*^2}\cong
{{\bf Z}}/2{{\bf Z}} \oplus
{{\bf Z}}/2{{\bf Z}}. \ $
The other part of the proof goes similarly as for  Theorem 3.

   \par \vskip 0.76cm

\begin{center}
{\Large IV. Criteria and Proof of Theorem 5}
\end{center}

  {\bf Proposition 4.1}. Let
  $C^\prime=C^\prime_{(-1)} :\;
 y^2=-x^4+(p+q)x^2-pq\ $ be as in Proposition 3.6.
 \par \vskip 0cm
  (i) $C^\prime $ has
  a rational integer point if and only if $p=3.$
\par \vskip 0cm

  (ii) $C^\prime $ has a rational point if and only if
  the following equation (I) or (II)  has a primary
   solution, i.e., a rational integer solution
  $(X,\ Y,\ S,\ T)$ with $(X, \ Y)=1$:
 $$ (I)\quad \left\{\begin{array}{l}
 X^2-pY^2 = S^2 \\
 X^2-qY^2 = -T^2 \end{array} \right. \qquad {{\rm or}}$$
 $$ (II)\quad \left\{\begin{array}{l}
 X^2-pY^2 = 2S^2 \\
 X^2-qY^2 = -2T^2 \end{array} \right. \hskip 1.2cm$$
 \hskip 0.9cm(iii) If
 $\ q=a^2+b^2,\; (a+\varepsilon)^2+(b+\delta)^2=c^2\ $ for
 some rational integer $\ a,\ b,\ c\ $ and
 $\ \varepsilon , \delta=\pm 1,\; $
 then the equation (I) in (ii)  has a primary
 solution, and $C^\prime $ has a rational point.

 \par \vskip 0.2cm

  {\bf Proof}. (i) If $p=3,\ q=5,\ $ then $(2,\ 1)$ is an
   integer point of $C^\prime \ $ (and obviously
   all the integer points of $C^\prime $ are just the
   four points $(\pm 2,\ \pm 1)$ )  .  On the other hand,
    if $(a,\ b)$ is a rational point of $C^\prime $,
    then $b^2=(a^2-p)(q-a^2)$, so
     $p<a^2<q=p+2,\ a^2=p+1, \ p=3.$

(ii) If $(a,\ b,\ c,\ d)$ is a primary solution of equation
(I) (or II), then
obviously $(a/b,\ cd/b^2)$ (or $(a/b,\ 2cd/b^2)$)
is a rational point of $C^\prime$.
 On the other hand,  if $(a/b,\ c/d)$ is a rational point of
 $C^\prime $ (and we may assume they are positive and
 $(a, b)=1$, $(c,d)=1$), then
  $(cb^2/d)^2=(a^2-pb^2)(qb^2-a^2)$.
 So  $u=cb^2/d\in {{\bf Z}},\ d|b^2 $,
 $\ (X,\ Y,\ Z) = (a,\ b,\ u)$    is a solution of
 $Z^2=(X^2-pY^2)(qY^2-X^2)$. It is easy to show
 $(a^2-pb^2,\ qb^2-a^2)=1$ or $2$.
 Put $u=u_1u_2$ or $2v_1v_2$
 with $(u_1,\ u_2)=1,\ (v_1,\ v_2)=1$, then
 $(a,\ b,\ u_1,\ u_2)$ or  $(a,\ b,\ v_1,\ v_2)$  is
 a primary solution of (I) or (II).
        \par \vskip 0cm

(iii)  Since $q>2$, we may assume $a+\varepsilon \neq 0$.
 Consider the equation
 $(a+\varepsilon)x^2 - 2(b+\delta )x - (a+\varepsilon)=0$,
 its discriminate is $4c^2$ so it has two rational roots
 $u,\ v$, so $\varepsilon u^2 - 2\delta u - \varepsilon
 = (1-u^2)a + 2ub$, $2(1+u^2)^2 =
 (\varepsilon u^2 + 2\delta u -\varepsilon)^2
 + (\varepsilon u^2 - 2\delta u - \epsilon)^2$,
 $q(1+u^2)^2 = ((1-u^2)a+2ub)^2+(2ua+(u^2-1)b)^2$.
 Thus $(2ua+(u^2-1)b,\ 1+u^2,\
 \varepsilon u^2+2\delta u-\varepsilon,\
 \varepsilon u^2-2\delta u -\varepsilon)$
 is a nontrivial rational
 solution of equation (I),
 from which a primary solution of (I)
 could be deduced easily.

\par \vskip 0.2cm
{\bf Proof of Theorem 5}. Since $p\equiv
 3\ ({{\rm mod}}\ 8)$,
so  by Proposition 3.6 we have
$-1\in S^{(\widehat{\varphi })}(E^\prime /{{\bf Q}})$. By
Proposition 4.1 (iii),
 $C^\prime $ has rational point, so
$-1\in E({{\bf Q}})/\widehat{\varphi}
(E^\prime ({{\bf Q}})).\ $ From the exact sequence
 \begin{center}
$0\rightarrow
{E({{\bf Q}})}/ {\widehat{\varphi}(E^\prime({{\bf Q}}))}
\rightarrow S^{(\widehat{\varphi  })}(E^\prime/{{\bf Q}})
\rightarrow $  $\Sh(E^\prime/{{\bf Q}} )[\widehat{\varphi}]
\rightarrow 0 ; $
     \end{center}
 we have
  $\Sh(E^\prime /{{\bf Q}} )[\widehat{\varphi }]=\{0\}.\;$
 And since
   $\Sh(E/{{\bf Q}} )[\varphi]=\{0\} ,\ $ so
   $\Sh(E/{{\bf Q}} )[2]=\{0\} .\ $
  Thus by Theorem 3 we know  $r=1$.

Now consider Examples in Section I.
Among all twin primes $p, q < 100,\ $
  only
$(3, 5)$ and $ (11, 13)$ satisfy the condition of
Theorem 5: $5=1^2+2^2,\ (1-1)^2+(2-1)^2=1^2;\;$
$13=2^2+3^2,\ (2+1)^2+(3+1)^2=5^2.\ $ So we know
${{\rm rank}} (E)=1.\ $
 For $(p,\ q)$=$(3,\ 5 )$, the equation (I) in
 Proposition 4.1(ii)
 has solution $(2,1,1,1)$.
 For $(p,\ q)$=$(11,\ 13 )$, the equation (I)
 has solution $(18, 5, 7, 1)$.

\vskip 1cm
\begin{center}{\bf {\large \bf R}eferences}\end{center}
 \baselineskip 0pt
\parskip 0pt
\begin{description}

\item[[1]] A. Bremner, On the equation $ y^2=x(x^2+p),\ $
in Number Theory and Applications
      (R. A. Mollin , ed.), Kluwer, Dordrecht, 3-23, 1989.
\item[[2]] A. Bremner and J. W. S. Cassels, On the equation
 $ y^2=x(x^2+p),\ $   Math. Comp. 42(1984), 257-264.
\item[[3]] R. J. Strocker and J. Top, On the equation
$ y^2=(x+p)(x^2+p^2),\ $ , Rocky Mountain J. of Math.
      24(1994), 1135-1161.
\item[[4]] J. P. Buhler, B. H. Gross and D. B. Zagier, On the
conjecture of Birch and Swinnerton-Dyer for an elliptic curve
of rank 3, Math. Comp. 44(1985), 473-481.
\item[[5]] J. H. Silverman, The Arithmetic of Elliptic Curves,
 GTM 106, Springer-Verlag, 1986.
\item[[6]] J. H. Silverman, The Advanced Topics in the
Arithmetic of Elliptic Curves, GTM 151, Springer-Verlag, 1994.
\item[[7]] J. Tate, Algorithm for determining the type of
singular fibber in an elliptic pencil, Modular
      Functions of One Variable IV, Lect. Notes Math. 476,
       Springer-Verlag, 1975, 33-52.
\item[[8]] F. Diamond, On deformation rings and Hecke rings,
 Annals of Math. 144(1996), 137-166.
\item[[9]] H. H. M$\ddot{u}$ller, H. Str$\ddot{o}$her and
H. G. Zimmer, Torsion subgroups of elliptic curves with
integral j-invariant over quadratic fields, J. Reine angew.
Math., 397(1989), 100-161.
\item[[10]] J. R. Merriman, S. SikSek and N. P. Smart,
Explicit 4-descents on an elliptic curves, Acta
      Arithmetica, LXXVII, 4(1996), 385-404.

\item[[11]] J. E. Cremona, Algorithms for modular elliptic
curves, Cambridge Univ. Press, 1992.
\item[[12]] J. W. S. Cassels, Lectures on Elliptic Curves,
LMS Student Texts, Cambridge Univ. Press, 1991.
\end{description}
\par  \vskip 0.3cm

{\Large T}SINGHUA {\Large U}NIVERITY \par \vskip 0.136cm
{\Large D}EPARTMENT OF
{\Large M}ATHEMATICAL {\Large S}CIENCES \par \vskip 0.136cm
{\Large B}EIJING 100084,  P. R. {\Large C}HINA \par
\vskip 0.16cm

E-mail address:\hskip 0.2cm  xianke@tsinghua.edu.cn

(This paper was published in:  Science In China, A45(2002.11),
No.11, 1372-1380)

\end{document}